# Undecidable proposition in PA and Diophantine equation

T. Mei

(Central China Normal University, Wuhan, Hubei PRO, People's Republic of China

E-Mail:　meitao@mail.ccnu.edu.cn　　meitaowh@public.wh.hb.cn )

**Abstract:** Based on the MRDP theorem concerning the Hilbert tenth problem, there is a corresponding Diophantine equation called *proof equation* for every formula of the First-order Peano Arithmetic (PA). A formula is provable in PA, if and only if the corresponding proof equation has solution. Based on proof equation, some famous sentences, e.g., the Gödel sentence, the Rosser sentence and the Henkin sentence, can be expressed by the form of Diophantine equation. It is proved that for every axiom and theorem in the PA, we can construct actually a corresponding Diophantine equation, for which we know that it has no solution, but this fact cannot be proved in the PA. This means that, for every axiom and theorem in the PA, we can construct actually a corresponding undecidable proposition. Finally, generalizing the idea of proof equation to any mathematical (set theoretical, number theoretical, algebraic, geometrical, topological, et al) proposition, a project translating the task seeking a proof of the mathematical proposition into solving a corresponding Diophantine equation is discussed.

## 1　The proof equation for any formula of PA

First-order Peano Arithmetic System (PA) is discussed in many standard reference books. In this paper, all the convention of the symbols, the terms, the formulas, the axiom schema, the deductive rules, the proof of a formula of PA and all the corresponding Gödel codes follows the Ref. [1]. Some systems weaker than PA, e.g., the system $Q$ [2, 3], are mentioned occasionally.

The solutions of an arbitrary Diophantine equation discussed in this paper are restricted to the set of natural numbers {0, 1, 2, …}. At first, we investigate simply the representability of Diophantine equation in PA.

According to the representability of recursive functions, an arbitrary polynomial $D(\{a_l\}, \{x_m\})$ with natural number coefficients can be represented by a formula $d(\{\boldsymbol{x}_{al}\}, \{\boldsymbol{x}_n\})$ in PA, where $\{a_l\}$, $\{x_m\}$ are the abbreviation of the parameters $a_1, a_2, \ldots, a_l$ and the variables $x_1, x_2, \ldots, x_m$, respectively. For obtaining $d(\{\boldsymbol{x}_{al}\}, \{\boldsymbol{x}_n\})$, what we need to do is only to translate every parameter $a_i$ in $D(\{a_l\}, \{x_m\})$ into $\boldsymbol{x}_{ai}$, $a_i^u$ into $\overbrace{\boldsymbol{x}_{ai} \times \boldsymbol{x}_{ai} \cdots \boldsymbol{x}_{ai}}^{u}$, every unknown $x_j$ in $D(\{a_l\}, \{x_m\})$ into $\boldsymbol{x}_j$, $x_j^v$ into $\overbrace{\boldsymbol{x}_j \times \boldsymbol{x}_j \cdots \boldsymbol{x}_j}^{v}$ and every constant $b_n$ into the term $\boldsymbol{0}''^{\cdots\prime}$ (There are $b_n$ successors "$\prime$"), which is abbreviated as $\boldsymbol{b}_n\prime$. The symbols in PA are expressed by boldface in this paper.

And, further, an arbitrary Diophantine equation $D(\{a_l\}, \{x_m\})=0$ can be rewritten as $D_L(\{a_l\}, \{x_m\})=D_R(\{a_l\}, \{x_m\})$, where both $D_L(\{a_l\}, \{x_m\})$ and $D_R(\{a_l\}, \{x_m\})$ are polynomials with natural number coefficients. Let the representations in PA of $D_L(\{a_l\}, \{x_m\})$ and $D_R(\{a_l\}, \{x_m\})$ be $d_L(\{x_{al}\}, \{x_n\})$ and $d_R(\{x_{al}\}, \{x_n\})$, respectively, then the representation of $D(\{a_l\}, \{x_m\})=0$ in PA is $d(\{x_{al}\}, \{x_n\})$, where $d(\{x_{al}\}, \{x_n\})$ is the formula $d_L(\{x_{al}\}, \{x_n\})=d_R(\{x_{al}\}, \{x_n\})$.

We therefore have

**Lemma 1.** For an arbitrary polynomial predicate $D(\{x_n\})=0$, There is a formula $d(\{x_n\})$ in PA such that

If $D(\{b_n\})=0$, then PA $\vdash d(\{b_n\}')$; If $D(\{b_n\})\neq 0$, then PA $\vdash \neg d(\{b_n\}')$.

**Corollary 1.** Let $d_{(1)}(\{x_l\})$ and $d_{(2)}(\{x_m\})$ are the representational formulas of the Diophantine equations $D_{(1)}(\{x_l\})=0$ and $D_{(2)}(\{x_m\})=0$, respectively, then representational formulas of the equations $D_{(1)}(\{x_{l+m}\})=(D_{(1)}(\{x_l\}))^2+(D_{(2)}(\{x_m\}))^2=0$ and $D_{(2)}(\{x_{l+m}\})=D_{(1)}(\{x_l\}) \cdot D_{(2)}(\{x_m\})=0$ are $d_{(1)}(\{x_{l0}\}')\wedge d_{(2)}(\{x_{m0}\}')$ and $d_{(1)}(\{x_{l0}\}')\vee d_{(2)}(\{x_{m0}\}')$, respectively.

**Lemma 2.** Let $d(\{x_n\})$ be the representational formula in PA of a Diophantine equation $D(\{x_n\})=0$, then

(1) If $D(\{x_n\})=0$ has solution then PA $\vdash \exists \{x_n\} d(\{x_n\})$;

(2) If PA is $\omega$-consistent, then if PA $\vdash \exists \{x_n\} d(\{x_n\})$, then the equation $D(\{x_n\})=0$ has solution.

*Proof*: (1) Immediate from Lemma 1.

(2) If PA $\vdash \exists \{x_n\} d(\{x_n\})$, namely, PA $\vdash \neg \forall \{x_n\} \neg d(\{x_n\})$. Suppose $D(\{x_n\})=0$ has no solution, from Lemma 1 we have PA $\vdash \neg d(\{l_n\}')$ for any $\{l_n\}'$. Hence, all the formulas $\neg d(\{l_n\}')$ for any $\{l_n\}'$ and $\neg \forall \{x_n\} \neg d(\{x_n\})$ are provable in PA, and hence PA is $\omega$-inconsistency.

**Remark 1.** (1) The notion of $\omega$-consistency is about a formula $A(x)$ of PA with just one free variable $x$. For the proof of Lemma 2, we should use the $n$-fold form of $\omega$-consistency called "$\omega^n$-consistency[3]", because there are many free variables $\{x_n\}$ but not one free variable in $d(\{x_n\})$. However, according to the following lemma in Ref.[3, §2]: "*Suppose $R_0\subseteq T$. If $T$ is $\omega$-consistency, the $T$ is $\omega^k$-consistency, for all $k$*", we can use only the condition of $\omega$-consistency in the proof of the result (2) of Lemma 2.

(2) For the result (2) of Lemma 2, the condition of $\omega$-consistency is necessary. In fact, there is a system in which all the formulas $\neg d(\{l_n\}')$ for any $\{l_n\}'$ and $\neg \forall \{x_n\} \neg d(\{x_n\})$ are provable. The system $D$ in the Ref.[3] can be as an example, because $D$ contains an infinite element. Just as Ref.[3] concludes, the system $D$ is $\omega$-inconsistency.

(3) If $D(\{x_n\})=0$ has no solution, from Lemma 1 we have PA $\vdash \neg d(\{l_n\}')$ for any $\{l_n\}'$. In this case, even if the condition of $\omega$-consistency is used, what we can obtain is only the conclusion "not PA $\vdash \neg \forall \{x_n\} \neg d(\{x_n\})$", e.g., "not PA $\vdash \exists \{x_n\} d(\{x_n\})$", we cannot obtain the conclusion "PA $\vdash \forall \{x_n\} \neg d(\{x_n\})$". In fact, there is equation $D(\{x_n\})=0$ which has no solution but the corresponding formula $\forall \{x_n\} \neg d(\{x_n\})$ is not provable in PA. We shall investigate such equations in next section.

**Theorem 1.** Defining a set $\text{Prf}_F$ for a formula $F$ of PA:

$\text{Prf}_F := \{a \in \text{Prf}_F | a \text{ is a Gödel code of a proof of the formula } F \text{ of PA}\}$,

then (1) The set $\text{Prf}_F$ is recursive.

(2) There is a Diophantine equation $D_{PE}(\#F, x_0, \{x_w\})=0$ such that $a \in \text{Prf}_F$ if and only if $D_{PE}(\#F, x_0, \{x_w\})=0$ has solution $x_0=a$ and $\{x_w\}=\{x_{w0}\}$, where $\#F$ is the Gödel code of the formula $F$. We call the equation $D_{PE}(\#F, x_0, \{x_w\})=0$ *the proof equation of the formula F in* PA.

(3) If $D_{PE}(\#F, x_0, \{x_w\})=0$ has solution, then PA $\vdash F$.

(4) There is a formula $d_{PE}((\#F)', \boldsymbol{x}_0, \{\boldsymbol{x}_w\})$ in PA as the representation of $D_{PE}(\#F, x_0, \{x_w\})=0$ such that

$$\text{If } D_{PE}(\#F, b_0, \{b_w\})=0, \text{ then PA} \vdash d_{PE}((\#F)', \boldsymbol{b}_0', \{\boldsymbol{b}_w\}');$$
$$\text{If } D_{PE}(\#F, b_0, \{b_w\})\neq 0, \text{ then PA} \vdash \neg\, d_{PE}((\#F)', \boldsymbol{b}_0', \{\boldsymbol{b}_w\}').$$

(5) If $D_{PE}(\#F, x_0, \{x_w\})=0$ has solution, then PA $\vdash \exists(\boldsymbol{x}_0, \{\boldsymbol{x}_n\})\, d_{PE}((\#F)', \boldsymbol{x}_0, \{\boldsymbol{x}_w\})$; If PA is $\omega$-consistent, then if PA $\vdash \exists(\boldsymbol{x}_0, \{\boldsymbol{x}_n\})\, d_{PE}((\#F)', \boldsymbol{x}_0, \{\boldsymbol{x}_w\})$, then the equation $D_{PE}(\#F, x_0, \{x_w\})=0$ has solution.

**Remark 2.** We are not going to give a formal proof since it is prolix. In fact, the thought of the proof is very simple.

At first, whatever how to implement a project of the Gödel code for PA, a basic character is that whether a given natural number $a$ is a Gödel code of a proof of a given formula $F$ is decidable, hence whether $a \in \text{Prf}_F$ is decidable, that means that $\text{Prf}_F$ is a recursive set.

And then, according to the MRDP theorem[4–6]: *A set is recursively enumerable if and only if it is Diophantine.* Hence, using the effective method given by the completely constructive process of proof of the MRDP theorem, we can find out a Diophantine equation $D_{PE}(\#F, a, \{x_w\})=0$ for the recursive set $\text{Prf}_F$ such that $a \in \text{Prf}_F$ if and only if $D_{PE}(\#F, a, \{x_w\})=0$ has at least solution $\{x_w\}=\{x_{w0}\}$ for $a$.

If we regard $a$ in the equation $D_{PE}(\#F, a, \{x_w\})=0$ as a unknown, then it is obvious that $D_{PE}(\#F, a, \{x_w\})=0$ has solution $\{x_w\}=\{x_{w0}\}$ if and only if the equation $D_{PE}(\#F, x_0, \{x_w\})=0$ has solution $x_0=a$ and $\{x_w\}=\{x_{w0}\}$.

(3) is obvious according to the definition of proof equation. And, further, if we obtain actually a group of obtained solution $x_0=a$ and $\{x_w\}=\{x_{w0}\}$ of the equation $D_{PE}(\#F, x_0, \{x_w\})=0$, then we obtain actually a proof of the proposition $F$ from the natural number $x_0=a$.

(4) and (5) of Theorem 1 are immediate from Lemma 1 and Lemma 2, respectively.

For an arbitrary formula $F$ of PA, starting from the corresponding proof equation $D_{PE}(\#F, x_0, \{x_w\})=0$ and the representational formula $d_{PE}((\#F)', \boldsymbol{x}_0, \{\boldsymbol{x}_w\})$ we can introduce a sequence of proof equations and the representational formulas about the formula $F$ as follows.

Using $F_{(0)}$ to denote $F$, $F_{(1)}$ to denote the formula $\exists(\boldsymbol{x}_0, \{\boldsymbol{x}_w\})\, d_{PE}((\#F)', \boldsymbol{x}_0, \{\boldsymbol{x}_w\})$, we have the corresponding $D_{PE}(\#F_{(1)}, x_{10}, \{x_{w1}\})=0$ and $d_{PE}((\#F_{(1)})', \boldsymbol{x}_{10}, \{\boldsymbol{x}_{w1}\})$. And, using $F_{(2)}$ to denote the formula $\exists(\boldsymbol{x}_{10}, \{\boldsymbol{x}_{1w}\})\, d_{PE}((\#F_{(1)})', \boldsymbol{x}_{10}, \{\boldsymbol{x}_{w2}\})$, we have the corresponding $D_{PE}(\#F_{(2)}, x_{20}, \{x_{w2}\})=0$ and $d_{PE}((\#F_{(2)})', \boldsymbol{x}_{20}, \{\boldsymbol{x}_{w2}\})$. Using $F_{(3)}$ to denote the formula $\exists(\boldsymbol{x}_{20}, \{\boldsymbol{x}_{w2}\})\, d_{PE}((\#F_{(2)})', \boldsymbol{x}_{20}, \{\boldsymbol{x}_{w2}\})$, we have the corresponding $D_{PE}(\#F_{(3)}, x_{30}, \{x_{w3}\})=0$ and $d_{PE}((\#F_{(3)})', \boldsymbol{x}_{30}, \{\boldsymbol{x}_{w3}\})$, ……, and so on. For these proof equations and the representational formulas, according to Theorem 1 we have

**Corollary 2.** If PA is $\omega$-consistent, then if an arbitrary equation $D_{PE}(\#F_{(k)}, x_{k0}, \{x_{wk}\})=0$ has solution, e.g., if PA $\vdash \exists(\boldsymbol{x}_{k0}, \{\boldsymbol{x}_{wk}\})\, d_{PE}(\#F_{(k)})', \boldsymbol{x}_{k0}, \{\boldsymbol{x}_{kw}\})$, then all equations $D_{PE}(\#F_{(n)}, x_{n0}, \{x_{wn}\})=0$ have solution and PA $\vdash F_{(n)}$, where $n=0, 1, 2, \dots k-1$; $F_{(0)} \equiv F$.

## 2  Some undecidable propositions expressed by Diophantine form in PA

At first, we give the following lemma without proof, since the proof is simple.

**Lemma 3.**

(1) All the following predicates are recursive,

① FC($a$) is true if and only if $a$ is the Gödel code of a formula of PA.

② Sub($a,b$) is true if and only if $a$ is the Gödel code of a formula $A(x)$ of PA with just one free variable $x$, namely, $a=\#A(x)$; $b$ is the Gödel code of a formula $A(a')$ obtained from $A(x)$ by substituting the numeral $a'$ for the free variable $x$.

③ Neg($a,b$) is true if and only if $a$ is the Gödel code of a formula $A(x)$ of PA, $b$ is the Gödel code of a formula $\neg A(x)$.

④ U($a,b$) is true if and only if $a$ is the Gödel code of a formula $A$ of PA, namely, $a=\#A$; $b$ is the Gödel code of a proof of the formula $A$.

⑤ V($a,b$) is true if and only if $a$ is the Gödel code of a formula $A$ of PA, $b$ is the Gödel code of a proof of the formula $\neg A$.

(2) Five Diophantine equations $D_{FC}(a, \{x_{w1}\})=0$, $D_{Sub}(a, b, \{x_{w2}\})=0$, $D_{Neg}(a, b, \{x_{w3}\})=0$, $D_U(a, b, \{x_u\})=0$ and $D_V(a, b, \{x_v\})=0$ can be constructed according to the MRDP theorem such that the predicates FC($a$), Sub($a, b$), Neg($a, b$), U($a, b$) and V($a, b$) are true if and only if the corresponding Diophantine equation $D_{FC}(a, \{x_{w1}\})=0$, $D_{Sub}(a, b, \{x_{w2}\})=0$, $D_{Neg}(a, b, \{x_{w3}\})=0$, $D_U(a, b, \{x_u\})=0$ and $D_V(a, b, \{x_v\})=0$ have solution, respectively.

(3) There are corresponding representational formulas $d_{FC}(x_0, \{x_{w1}\})$, $d_{Sub}(x_1, x_2, \{x_{w2}\})$, $d_{Neg}(x_1, x_2, \{x_{w3}\})$, $d_U(x_1, x_2, \{x_u\})$ and $d_V(x_1, x_2, \{x_v\})$ in PA for the equations $D_{FC}(a, \{x_{w1}\})=0$, $D_{Sub}(a, b, \{x_{w2}\})=0$, $D_{Neg}(a, b, \{x_{w3}\})=0$, $D_U(a, b, \{x_u\})=0$ and $D_V(a, b, \{x_v\})=0$, respectively.

In fact, $D_U(a, b, \{x_u\})=(D_{FC}(a, \{x_{w1}\}))^2+(D_{PE}(a, b, \{x_w\}))^2$, $D_V(a, b, \{x_v\})=(D_{FC}(a, \{x_{w1}\}))^2+(D_{Neg}(a, x_0, \{x_{w3}\}))^2+(D_{PE}(x_0, b, \{x_w\}))^2$, both $d_U(x_1, x_2, \{x_u\})$ and $d_V(x_1, x_2, \{x_v\})$ thus can be expressed by $d_{FC}(x_0, \{x_{w1}\})$, $d_{Sub}(x_1, x_2, \{x_{w2}\})$ and $d_{Neg}(x_1, x_2, \{x_{w3}\})$ according to Corollary 1.

Now we consider a formula $g_1(x)$ of PA with just one free variable $x$:

$$g_1(x) \leftrightarrow \forall (x_1, \{x_{w2}\}) (d_{Sub}(x, x_1, \{x_{w2}\}) \rightarrow \neg \exists (x_2, \{x_w\}) d_{PE}(x_1, x_2, \{x_w\})),$$

substituting the numeral $(\#g_1(x))'$ of the Gödel code of the formula $g_1(x)$ for the free variable $x$ in $g_1(x)$, we obtain a formula $G_1$:

$$G_1 \leftrightarrow \forall (x_1, \{x_{w2}\}) (d_{Sub}((\#g_1(x))', x_1, \{x_{w2}\}) \rightarrow \neg \exists (x_2, \{x_w\}) d_{PE}(x_1, x_2, \{x_w\})).$$

$G_1$ has just the meaning of the Gödel sentence "I am not provable in PA", for which we have:

**Theorem 2** If PA is consistent, then $G_1$ is not provable in PA; if PA is $\omega$-consistent, then $\neg G_1$ is not provable in PA.

*Proof*: At first, notice that Sub($\#g_1(x),\#G_1$) is true, namely, the equation $D_{Sub}(\#g_1(x), \#G_1, \{x_{w2}\})=0$ has solution, where $\#G_1$ is the Gödel code of $G_1$.

Now we prove that $G_1$ is not provable in PA. Contrarily, suppose PA $\vdash G_1$, then the equation $D_{PE}(\#G_1, x_1, \{x_w\})=0$ has solution, and then the equation $D_{G1}=(D_{Sub}(\#g_1(x), \#G_1, \{x_{w2}\}))^2$ $+(D_{PE}(\#G_1, x_2, \{x_w\}))^2=0$ has solution; Assuming a group of solution are $\{x_{w2}\}=\{x_{w20}\}$, $x_2=x_{20}, \{x_w\}=\{x_{w0}\}$, according to Theorem 1(3), PA $\vdash$ $d_{Sub}((\#g_1(x))'$, $(\#G_1)'$, $\{x_{w20}\}')\wedge d_{PE}((\#G_1)'$, $(x_{20})'$, $\{x_{w0}\}')$, hence PA $\vdash$ $\exists (x_1, \{x_{w2}\}, x_2, \{x_w\})$ $(d_{Sub}((\#g_1(x))'$, $x_1, \{x_{w2}\})\wedge d_{PE}(x_1, x_2, \{x_w\}))$. This is a contradiction with PA $\vdash$ $G_1$, because $G_1\leftrightarrow\neg\exists(x_1, \{x_{w2}\}, x_2, \{x_w\})$ $(d_{Sub}((\#g_1(x))'$, $x_1, \{x_{w2}\})\wedge d_{PE}(x_1, x_2, \{x_w\}))$.

Now we prove that $\neg G_1$ is not provable in PA. Contrarily, suppose PA $\vdash \neg G_1$, namely, PA $\vdash \exists(x_1, \{x_{w2}\}, x_2, \{x_w\})$ $(d_{Sub}((\#g_1(x))'$, $x_1, \{x_{w2}\})\wedge d_{PE}(x_1, x_2, \{x_w\}))$, according to Theorem 1(3) and Theorem 1 (5), if PA is $\omega$-consistent, then the equation $D_{G1}=(D_{Sub}(\#g_1(x), \#G_1, \{x_{w2}\}))^2 +(D_{PE}(\#G_1, x_2, \{x_w\}))^2=0$ has solution. However, we have already proved that equation $D_{G1}$ has no solution, a contradiction thus occurs. The Theorem 2 is proved.

**Corollary 3**. If PA is consistent, then the equation $D_{G1}=(D_{Sub}(\#g_1(x), \#G_1, \{x_{w2}\}))^2+(D_{PE}(\#G_1, x_2, \{x_w\}))^2=0$ has no solution, but the corresponding representational formula $\forall (x_1, \{x_{w2}\}, x_2, \{x_w\})$ $\neg$ $(d_{Sub}((\#g_1(x))'$, $x_1, \{x_{w2}\})\wedge d_{PE}(x_1, x_2, \{x_w\}))$ is not provable in PA.

*Proof*: Immediate from Theorem 2.

**Lemma 4 (Diagonal Theorem)**[2]. For an arbitrary formula $α(x)$ of PA with just one free variable $x$, there is a sentence $A$ such that PA $\vdash A\leftrightarrow α(\#A)$.

Now consider a formula $g_2(x)$ of PA with just one free variable $x$:

$g_2(x)\leftrightarrow \forall \{x_{w1}\}$ $(d_{FC}(x, \{x_{w1}\})\rightarrow \neg \exists(x_1, \{x_w\})$ $d_{PE}(x, x_1, \{x_w\}))$ ;

According to the Lemma 4, there is a sentence $G_2$ such that

PA $\vdash G_2\leftrightarrow \forall \{x_{w1}\}$ $(d_{FC}((\#G_2)'$, $\{x_{w1}\})\rightarrow \neg \exists(x_1, \{x_w\})$ $d_{PE}((\#G_2)'$, $x_1, \{x_w\}))$ .

$G_2$ has also just the meaning of the Gödel sentence "I am not provable in PA", for which we have:

**Theorem 3.** If PA is consistent, then $G_2$ is not provable in PA; if PA is $\omega$-consistent, then $\neg G_2$ is not provable in PA.

*Proof*: Notice that $FC(\#G_2)$ is true, namely, the equation $D_{FC}(\#G_2, \{x_{w1}\})=0$ has solution. The rest process of the proof is similar to that of the Theorem 2.

**Corollary 4.** If PA is consistent, then the equation $D_{G2}=(D_{FC}(\#G_2, \{x_{w1}\}))^2+(D_{PE}(\#G_2, x_1, \{x_w\}))^2=0$ has no solution, but the corresponding representational formula $\forall (\{x_{w1}\}, x_1, \{x_w\})$ $\neg$ $(d_{FC}((\#G_2)'$, $\{x_{w1}\})\wedge d_{PE}(\#G_2, x_1, \{x_w\}))$ is not provable in PA.

The proof of the Corollary 4 is similar to that of the Corollary 3.

**Corollary 5.** If PA is consistent, then the equation $D_{PE}(\#G_2, x_1, \{x_w\})=0$ has no solution, but the corresponding representational formula $\forall (x_1, \{x_w\})$ $\neg$ $d_{PE}((\#G_2)'$, $x_1, \{x_w\})$ is not provable in PA.

*Proof*: If PA $\vdash \forall (x_1, \{x_w\})$ $\neg$ $d_{PE}((\#G_2)'$, $x_1, \{x_w\})$, then PA $\vdash (\forall (x_1, \{x_w\})$ $\neg$ $d_{PE}((\#G_2)'$, $x_1, \{x_w\}))\vee Q$ for an arbitrary formula $Q$, and in particular

①      PA $\vdash (\forall (x_1, \{x_w\})$ $\neg$ $d_{PE}((\#G_2)'$, $x_0, \{x_w\}))\vee \neg \exists \{x_{w1}\}$ $d_{FC}((\#G_2)'$, $\{x_{w1}\})$ .

On the other hand, according to the theorem "$\forall x(α\rightarrow β)\leftrightarrow \exists x α\rightarrow β$ ($x$ does not occur freely in $β$)" in PA, we have

②      $G_2\leftrightarrow (\forall (x_1, \{x_w\})$ $\neg$ $d_{PE}((\#G_2)'$, $x_0, \{x_w\}))\vee \neg \exists \{x_{w1}\}$ $d_{FC}((\#G_2)'$, $\{x_{w1}\})$ .

From ① and ②, we have PA $\vdash G_2$, this is a contradiction with Theorem 3.

Now we consider the formula $r(x)$:

$r(x) \leftrightarrow \forall (x_1, \{x_u\}) (d_U(x, x_1, \{x_u\}) \to \exists (x_2, \{x_v\}) (x_2+\{x_v\} \leqslant x_1+\{x_u\} \wedge d_V(x, x_2, \{x_v\})))$

where $x_1+\{x_u\}$ and $x_2+\{x_v\}$ are the abbreviation of $x_1+\Sigma x_u$ and $x_2+\Sigma x_v$, respectively. According to the Lemma 4, there is a sentence $R$ such that

PA ⊢ $R \leftrightarrow \forall (x_1, \{x_u\}) (d_U((\#R)', x_1, \{x_u\}) \to$
$\exists (x_2, \{x_v\}) (x_2+\{x_v\} \leqslant x_1+\{x_u\} \wedge d_V((\#R)', x_2, \{x_v\})))$.

The meaning of $R$ is similar to the Rosser sentence "If I am provable in PA, then for every proof of me, there is a shorter proof of my negation", for which we have:

**Theorem 4.** If PA is consistent, then both $R$ and $\neg R$ are not provable in PA.

*Proof*: The proof is similar to the Rosser's proof[1,2,7]. At first, we prove that $R$ is not provable in PA. Contrarily, suppose PA ⊢ $R$, then the equation $D_U(\#R, x_1, \{x_u\})=0$ has solution, assuming a group of solution are $x_1=x_{10}, \{x_u\}=\{x_{u0}\}$, then PA ⊢ $d_U((\#R)', x_{10}', \{x_{u0}\}')$. According to the axiom $\forall x\, a(x) \to a(x/t)$ and using the MP rule, we obtain:

$\quad$ PA ⊢ $\exists (x_2, \{x_v\}) (x_2+\{x_v\} \leqslant x_{10}'+\{x_{u0}\}' \wedge d_V((\#R)', x_2, \{x_v\}))$ $\quad$ ①

On the other hand, if PA is consistent, then now $\neg R$ is not provable, and then the equation $D_V(\#R, x_2, \{x_v\})=0$ has no solution, hence PA ⊢ $\neg d_V((\#R)', x_{20}', \{x_{v0}\}')$ for an arbitrary numeral $x_{20}'$ and $\{x_{v0}\}'$; In particular, PA ⊢ $\neg d_V((\#R)', 0, \{0\})$, PA ⊢ $\neg d_V((\#R)', 0', \{0\})$, …, PA ⊢ $\neg d_V((\#R)', x_{20}', \{x_{v0}\}')$, where $x_{20}'+\{x_{v0}\}' \leqslant x_{10}'+\{x_{u0}\}'$; and then PA ⊢ $\neg d_V((\#R)', 0, \{0\}) \wedge \neg d_V((\#R)', 0', \{0\}) \wedge … \wedge \neg d_V((\#R)', x_{20}', \{x_{v0}\}')$, where $x_{20}'+\{x_{v0}\}' \leqslant x_{10}'+\{x_{u0}\}'$; and then PA ⊢ $\forall x_2, \{x_v\} (x_2+\{x_v\} \leqslant x_{10}'+\{x_{u0}\}' \to \neg d_V((\#R)', x_2, \{x_v\}))$, namely,

$\quad$ PA ⊢ $\neg \exists (x_2, \{x_v\}) (x_2+\{x_v\} \leqslant x_{10}'+\{x_{u0}\}' \wedge d_V((\#R)' x, x_2, \{x_v\}))$

This is a contradiction with ①.

Now we proof that $\neg R$ is not provable in PA. Contrarily, suppose PA ⊢ $\neg R$, then the equation $D_V(\#R, x_2, \{x_v\})=0$ has solution, assuming a group of solution are $x_2=x_{20}, \{x_v\}=\{x_{v0}\}$, then

(1) $\quad$ PA ⊢ $d_V((\#R)', x_{20}', \{x_{v0}\}')$

(2) $\quad$ PA $\cup \{x_{20}'+\{x_{v0}\}' \leqslant x_1+\{x_u\}\}$ ⊢ $\{x_{20}'+\{x_{u0}\}' \leqslant x_1+\{x_u\} \wedge d_V((\#R)', x_{20}', \{x_{v0}\}')\}$

(3) $\quad$ PA $\cup \{x_{20}'+\{x_{v0}\}' \leqslant x_1+\{x_u\}\}$ ⊢ $\exists x_2, \{x_v\}(x_2+\{x_v\} \leqslant x_1+\{x_u\} \wedge d_V((\#R)', x_2, \{x_v\}))$

(4) $\quad$ PA ⊢ $x_{20}'+\{x_{v0}\}' \leqslant x_1+\{x_u\} \to \exists x_2, \{x_v\}(x_2+\{x_v\} \leqslant x_1+\{x_u\} \wedge d_V((\#R)', x_2, \{x_v\}))$

On the other hand, if PA is consistent, then now $R$ is not provable, and then the equation $D_U(\#R, x_1, \{x_u\})=0$ has no solution, hence PA ⊢ $\neg d_U((\#R)', x_{10}', \{x_{u0}\}')$ for an arbitrary numeral $x_{10}'$ and $\{x_{u0}\}'$; In particular, PA ⊢ $\neg d_U((\#R)', 0, \{0\})$, PA ⊢ $\neg d_U((\#R)', 0', \{0\})$, …, PA ⊢ $\neg d_U((\#R)', x_{10}', \{x_{u0}\}')$, where $x_{10}'+\{x_{u0}\}' \leqslant x_{20}'+\{x_{v0}\}'$; and then PA ⊢ $\neg d_U((\#R)', 0, \{0\}) \wedge \neg (d_U((\#R)', 0', \{0\}) \wedge … \wedge \neg d_U((\#R)', x_{10}', \{x_{u0}\}')$, where $x_{10}'+\{x_{u0}\}' \leqslant x_{20}'+\{x_{v0}\}'$; hence

(5) $\quad$ PA ⊢ $x_1+\{x_u\} \leqslant x_{20}'+\{x_{v0}\}' \to \neg d_U((\#R)', x_1, \{x_u\})$

(6) $\quad$ PA ⊢ $d_U((\#R)', x_1, \{x_u\}) \to \neg (x_1+\{x_u\} \leqslant x_{20}'+\{x_{v0}\}')$

(7) $\quad$ PA ⊢ $\neg (x_1+\{x_u\} \leqslant x_{20}'+\{x_{v0}\}') \to x_1+\{x_u\} > x_{20}'+\{x_{v0}\}'$

(8) $\quad$ PA ⊢ $x_1+\{x_u\} > x_{20}'+\{x_{v0}\}' \to x_1+\{x_u\} \geqslant x_{20}'+\{x_{v0}\}'$

From (6), (7), (8) and (4) we have

(9) $\quad$ PA ⊢ $d_U((\#R)', x_1, \{x_u\}) \to \exists (x_2, \{x_v\}) (x_2+\{x_v\} \leqslant x_1+\{x_u\} \wedge d_V((\#R)', x_2, \{x_v\}))$

According to the Gen rule $a$ ⊢ $\forall x\, a$, we obtain

PA ⊢ ∀ ($x_1$,{$x_u$}) ($d_U$((#R)′, $x_1$, {$x_u$})→∃ ($x_2$,{$x_v$}) ($x_2$+{$x_v$}≤$x_1$+{$x_u$}∧$d_V$((#R)′, $x_2$, {$x_v$})))

That means PA ⊢ $R$ and is a contradiction with PA ⊢ ¬ $R$.

**Corollary 6.** If PA is consistent, then the equation $D_U$(#R, $x_1$, {$x_u$})=($D_{FC}$(#R, {$x_{w1}$}))$^2$ +($D_{PE}$(#R, $x_1$, {$x_w$}))$^2$=0 has no solution, but the corresponding representational formula ∀ ($x_1$,{$x_u$}) ¬ $d_U$((#R)′, $x_1$, {$x_u$}) is not provable in PA.

*Proof*: The conclusion "The equation $D_U$(#R, $x_1$, {$x_u$})=0 has no solution" has proved in the proof of Theorem 4, if PA ⊢ ∀ ($x_1$,{$x_u$}) ¬ $d_U$((#R)′, $x_1$, {$x_u$}), then PA ⊢ ∀ ($x_1$,{$x_u$}) (¬ $d_U$((#R)′, $x_1$, {$x_u$})∨$Q$) for an arbitrary formula $Q$, and in particular PA ⊢ ∀ ($x_1$,{$x_u$}) (¬ $d_U$((#R)′, $x_1$, {$x_u$})∨∃ ($x_2$,{$x_v$}) ($x_2$+{$x_v$}≤$x_1$+{$x_u$}∧$d_V$($x$, $x_2$, {$x_v$}))), this means PA ⊢ $R$ and is a contradiction with Theorem 4.

**Remark 3.** For obtaining the conclusion "There are some Diophantine equations for which we know that they have no solutions, but these facts cannot be proved in the PA", from the proof of Corollary 3, 4, 5 and 6 we see that what we use is only "PA is consistent" but not "PA is ω-consistent". This is not inconsistent with the conclusions in the Ref. [3]. We state the Theorem 2 in the Ref. [3] as follows: "*There is a finitely axiomatizable consistent extension D of Q in which all Diophantine sentence are decidable.*" However, the system $D$ is in fact essentially different from $Q$ or PA, because, from the axiom $D_3$ (∃$z$∀$x$ ($x$+$z$=$z$)) in the Ref. [3] we see, $D$ is a system containing an infinite element. In such system, "all Diophantine sentence are decidable" is not surprising; just as the discussion in the Ref. [3], in a nonstandard model $M$ (which is used as the model of $D$) containing an infinite element ∞, an arbitrary Diophantine equation $P$({$x_m$})=$Q$({$x_m$}) have a trivial solution ∞, unless $P$ or $Q$ are polynomials of degree zero. And, further, even if the nonstandard model $M$ can be used as the model of $Q$ or PA, and hence all provable formulas are true in $M$, we still have not the conclusion that all true sentences are provable in $Q$ or PA (This is just what the incompleteness theorem shows). In fact, the true sentence "all Diophantine sentence are decidable" in $M$ is surely not a theorem in the system $Q$ or PA, because there is not an infinite element in $Q$ or PA, namely, there is not an axiom like $D_3$ (∃$z$∀$x$ ($x$+$z$=$z$)) to guarantee the existence of an infinite element in $Q$ or PA.

**Lemma 5.** Let $B(x)$ be ∃($x_1$, {$x_w$}) $d_{PE}$($x$, $x_1$, {$x_w$}), then if PA is consistent, then for an arbitrary formula $A$,

(P1)  If PA ⊢ $A$, then PA ⊢ $B$(#$A$);

(P2)  PA ⊢ $B$(#$A$) →$B$(#($B$(#$A$))).

If PA is ω-consistent, then for arbitrary formulas $A_1$, $A_2$,

(P3)  PA ⊢ $B$(#($A_1$→ $A_2$)) →($B$(#$A_1$)→ $B$(#$A_2$)).

*Proof*: (P1) Immediate from the definitions of proof equation and $B(x)$ and Theorem 1(5).

(P2) If PA ⊢ $B$(#$A$), then the equation $D_{PE}$(#($B$(#$A$)), $x_1$, {$x_w$})=0 has solution, assuming a group of solution are $x_1$=$x_{10}$,{$x_w$}={$x_{w0}$}, then PA ⊢ $d_{PE}$((#($B$(#$A$)))′, $x_{10}$′, {$x_{w0}$}′), and then PA ⊢ ∃$x_1$, {$x_w$} $d_{PE}$((#($B$(#$A$)))′, $x_1$, {$x_w$}), namely, PA ⊢ $B$(#($B$(#$A$))).

(P3) If PA ⊢ $B$(#($A_1$ → $A_2$)), according to Theorem 1(5), if PA is ω-consistent, then the equation $D_{PE}(\#(A_1 \to A_2), x_1, \{x_w\})=0$ has solution, according to Theorem 1(3) we have PA ⊢ $A_1 \to A_2$. Similarly, if PA ⊢ $B(\#A_1)$, then PA ⊢ $A_1$. Using the MP rule we have PA ⊢ $A_2$, then we have PA ⊢ $B(\# A_2)$. The (P3) of Lemma 5 is proved.

**Lemma 6 (Löb's Theorem)**[2, 8]. Let $B(x)$ be $\exists(x_1, \{x_w\}) \, d_{PE}(x, x_1, \{x_w\})$, then for an arbitrary formula $A$, if PA ⊢ $B(\#A) \to A$, the PA ⊢$A$.

**Remark 4.** The properties (P1) ~ (P3) in Lemma 5 is so called **"Löb's derivability conditions"**[2]. According to Lemma 5, both Löb's Theorem and the provability of the Henkin sentence[2, 8, 9] $H \leftrightarrow \exists(x_1, \{x_w\}) \, d_{PE}((\#H)', x_1, \{x_w\})$ in PA can be proved easily. However, in literatures, what used condition in the proof that the provability predicate $B(x)$ satisfies (P1) ~ (P3) is only "PA is consistent", but we have to use the condition "PA is ω-consistent" for the proof of the property (P3) in Lemma 5, this is also the reason that we use the condition "PA is ω-consistent" in the both proofs of below Lemma 7 and Theorem 5.

**Lemma 7.** If PA is ω-consistent and $B(x) \leftrightarrow \exists(x_1, \{x_w\}) \, d_{PE}(x, x_1, \{x_w\})$, then for an arbitrary formula $A$, if $A$ is an axiom or a theorem in PA, then not PA ⊢ ¬ $B(\#(¬ A))$.

*Proof*: Suppose PA ⊢ ¬ $B(\#(¬ A))$, then PA ⊢ $B(\#(¬ A)) \to F$ for any sentence $F$, and in particular PA ⊢ $B(\#(¬ A)) \to ¬ A$, and according to Lemma 6 we have PA ⊢ ¬ $A$. On the other hand, $A$ is an axiom or a theorem in PA, we therefore have also PA ⊢ $A$. PA is thus inconsistent.

In fact, Lemma 7 is the abstract form[2] of the Gödel's Second Incompleteness Theorem. For the below proof of Theorem 5, what we need is only this abstract form.

**Theorem 5.** If PA is ω-consistent, then for any axiom or theorem in the PA, there is a corresponding Diophantine equation, for which we know that it has no solutions, but this fact cannot be proved in the PA.

*Proof*: For any axiom or theorem $A$, ¬ $A$ is not provable in PA, namely, the equation $D_{PE}(\#(¬ A), x_1, \{x_w\})=0$ has no solution; however, if PA ⊢ $\forall (x_1, \{x_w\}) \, ¬ \, d_{PE}((\#(¬ A))', x_1, \{x_w\})$, namely, PA ⊢ ¬ $B(\#(¬ A))$, this is a contradiction with Theorem 5.

Notice that Theorem 5 means that, for every axiom and theorem in the PA, there is a corresponding undecidable proposition. And, further, the corresponding undecidable proposition can be constructed actually and has the Diophantine form. Theorem 5 shows that undecidable propositions in PA are at least as many as the axioms and theorems of PA.

## 3  A project translating the task seeking a proof of the mathematical proposition into solving a corresponding Diophantine equation

Generalizing the idea of proof equation to any mathematical (set theoretical, number theoretical, algebraic, geometrical, topological, et al) proposition, in principle, the task seeking a proof of the mathematical proposition may be translated into solving a corresponding Diophantine equation. Of course, maybe it is also difficult to solve a Diophantine equation, however, this is another kind of question.

The method of constructing the corresponding Diophantine equation for solving a mathematical problem has been discussed, see, for example, Ref. [3, 4, 5, 10]. However, our project differs from the method given by the Ref. [3, 4, 5, 10]. In the Ref. [3, 4, 5, 10], special analysis and discussion are necessary for every special mathematical problem, respectively. What project we present is universal.

As the first step of the project, for a proposition $F$ of a mathematical system $S$, we construct a corresponding formal system **S**, namely, S is as a model of **S**. Especially, there is a formula **F** in **S** to express the proposition $F$ of $S$.

Notice that what we want to do is to seek a proof for a proposition but not for studying a formal system, hence, the construction of **S** is comparatively arbitrary. For example, if the formula **F** is not provable in **S**, notice that it does not mean that the proposition $F$ is also not provable in the mathematical system $S$, then we can add new symbol, new axiom, extend to second-order system from first-order system, etc, to obtain another "stronger" formal system **S′** than **S**, and try to seek a proof for the proposition $F$ using the "stronger" system **S′**.

And, further, we can add some known theorems into **S** as axioms, because (1) If we start from the axioms to prove a formula of **S**, then maybe the proof is too lengthy; (2) Mathematicians start usually from known theorems to prove a mathematical proposition.

Of course, although the construction of the formal system **S** is comparatively arbitrary, some basic conditions must be satisfied. What we ask is merely whether a given symbol is one of the symbols, a sequence consisting of the symbols of **S** is a formula, and a sequence consisting of the formulas of **S** is a proof of a formula of **S** are decidable after the symbols, formulas and proof of a given formula in **S** are defined.

And then, we implement the Gödel code for **S**, for which what we ask is merely that whether a given natural number is the Gödel code of a symbol, a given natural number is the Gödel code of a formula, and a given natural number is the Gödel code of a proof of a formula of **S** are decidable.

For the formula **F** of **S**, we define a set $Prf_F$ consisting of natural numbers such that $a \in Prf_F$ if and only if $a$ is a Gödel code of a proof of **F**.

Because there is an effective method for determining whether a given natural number is the Gödel code of a proof of a given formula of **S**, the set $Prf_F$ is thus recursive. And then, using the effective method given by the completely constructive process of proof of the MRDP theorem, we can find out a Diophantine equation $D_{PE}(\#\boldsymbol{F}, a, \{x_w\})=0$ for the recursive set $\varXi_F$ such that $a \in Prf_F$ if and only if $D_{PE}(\#\boldsymbol{F}, a, \{x_w\})=0$ has at least a solution $\{x_w\}=\{x_{w0}\}$ for $a$.

The last step of the project is that we regard $a$ in the equation $D_{PE}(\#\boldsymbol{F}, a, \{x_w\})=0$ as a unknown variable and try to solve the equation $D_{PE}(\#\boldsymbol{F}, x_0, \{x_w\})=0$. We call the equation $D_{PE}(\#\boldsymbol{F}, x_0, \{x_w\})=0$ *the proof equation of the formula **F** in **S***.

There are great differences between equations $D_{PE}(\#\boldsymbol{F}, a, \{x_w\})=0$ and $D_{PE}(\#\boldsymbol{F}, x_0, \{x_w\})=0$. For example, whether the equation $D_{PE}(\#\boldsymbol{F}, a, \{x_w\})=0$ has solution for $a$ is decidable, because the set $Prf_F$ is recursive. And according to the algorithm deciding whether a given natural number belongs to the set $Prf_F$ if we know $a \in Prf_F$ then we also know that $D_{PE}(\#\boldsymbol{F}, a, \{x_w\})=0$ has solution for $a$; Conversely, if we know $a \notin Prf_F$ then we also know the solution of the equation $D_{PE}(\#F, a, \{x_w\})=0$ does not exist. However, if we want to use the same method to decide whether there is a solution of the equation $D_{PE}(\#F, x_0, \{x_w\})=0$, then if it has no solution then we have to check all natural numbers for the variable $x_0$ and the process never stop, we therefore do not know actually whether the solution exists.

For the proof equation $D_{PE}(\#F, x_0, \{x_w\})=0$, we have the following conclusions:

(1) $D_{PE}(\#F, x_0, \{x_w\})=0$ has a group of solution $x_0=a$, $\{x_w\}=\{x_{w0}\}$ if and only if $D_{PE}(\#F, a, \{x_w\})=0$ has a group of solution $\{x_w\}=\{x_{w0}\}$.

(2) If $D_{PE}(\#F, x_0, \{x_w\})=0$ has solution, then the formula $F$ is provable in $S$ and a proof of $F$ is actually obtained from the natural number $x_0=a$, we therefore obtain a proof of the proposition $F$ of the mathematical system $S$ and $F$ is turned into a theorem of $S$.

(3) If there is not any solution of the equation $D_{PE}(\#F, x_0, \{x_w\})=0$, namely, the formula $F$ is not provable in $S$, then we can further consider the formula $\neg F$ in $S$ and corresponding Diophantine equation $D_{PE}(\#(\neg F), x_0, \{x_w\})=0$, where $\#(\neg F)$ is the Gödel code of the formula $\neg F$.

① If $D_{PE}(\#(\neg F), x_0, \{x_w\})=0$ has solution, then the formula $\neg F$ is provable in $S$ and a proof of $\neg F$ is actually obtained from the natural number $x_0=a$, we therefore obtain a proof of the proposition $\neg F$ of the mathematical system $S$ and $\neg F$ is turned into a theorem of $S$.

② If there is also not any solution of the equation $D_{PE}(x_0, \#(\neg F), \{x_w\})=0$, then the formula $\neg F$ is also not provable in $S$. Now that both formulas $F$ and $\neg F$ are not provable in $S$, the formula $F$ is undecidable in $S$.

In particular, if a formula $F$ with the universal quantifier $\forall$ is undecidable in $S$, then we know actually sometimes that the corresponding proposition $F$ is true.

(4) If $F$ is undecidable in $S$, notice that it does not mean that the proposition $F$ is also not provable in the mathematical system $S$, then we have to try to construct another "stronger" formal system $S'$ than $S$ and repeat the above steps if we want to try to seek a proof to the proposition $F$ via the method of solving proof equation.